\documentclass{iopart}  

\eqnobysec

\usepackage{xypic}

\usepackage{amsthm}
\theoremstyle{plain}
\newtheorem{proposition}{Proposition}[section]
\newtheorem{corollary}[proposition]{Corollary}
\newtheorem{lemma}[proposition]{Lemma}
\newtheorem{theorem}[proposition]{Theorem}

\theoremstyle{definition}
\newtheorem{definition}[proposition]{Definition}

\newtheorem{remark}[proposition]{Remark}

\newcommand{\R}{\mathbf{R}}
\newcommand{\abs}[1]{\left\vert #1 \right\vert}

\newcommand{\sgn}{\mathop{\mathrm{sgn}}\nolimits}
\newcommand{\adj}{\mathop{\mathrm{adj}}\nolimits}
\newcommand{\degree}{\mathop{\mathrm{deg}}\nolimits}
\newcommand{\diag}{\mathop{\mathrm{diag}}\nolimits}

\newcommand{\pipos}{\Pi_{\ge0}}
\newcommand{\pineg}{\Pi_{<0}}
\newcommand{\pil}{\Pi^l}

\newcommand{\matrixM}{\mathcal{M}}
\newcommand{\matrixL}{\mathcal{L}}
\newcommand{\matrixX}{\mathcal{X}}
\newcommand{\matrixT}{\mathcal{T}}

\renewenvironment{pmatrix}{\left(\begin{array}{cccccccccc}}{\end{array}\right)}
\newenvironment{vmatrix}{\left\vert\begin{array}{cccccccccc}}{\end{array}\right\vert}
\newcommand{\binom}[2]{{#1\choose #2}}
\newcommand{\iint}{\int\!\!\!\int}

\begin{document}

\title{The inverse spectral problem for the discrete cubic string}

\author{Jennifer Kohlenberg\dag{}, Hans Lundmark\ddag{} and Jacek Szmigielski\S}
\address{\dag
  Institute of Biomaterials and Biomedical Engineering,
  University of Toronto,
  Toronto, ON, M5S 3G9, Canada
}
\address{\ddag
  Department of Mathematics,
  Link{\"o}ping University,
  SE-581 83 Link{\"o}ping, Sweden}
\address{\S
  Department of Mathematics and Statistics,
  University of Saskatchewan,
  106 Wiggins Road, Saskatoon, Saskatchewan, S7N 5E6, Canada}
\eads{\mailto{jennifer.kohlenberg@utoronto.ca}, \mailto{halun@mai.liu.se}, \mailto{szmigiel@math.usask.ca}}

\begin{abstract}
  Given a measure $m$ on the real line or a finite interval, the
  \emph{cubic string} is the third order ODE $-\phi'''=zm\phi$ where
  $z$ is a spectral parameter. If equipped with Dirichlet-like
  boundary conditions this is a nonselfadjoint boundary value problem
  which has recently been shown to have a connection to the
  Degasperis-Procesi nonlinear water wave equation. In this paper we
  study the spectral and inverse spectral problem for the case of
  Neumann-like boundary conditions which appear in a high-frequency
  limit of the Degasperis--Procesi equation. We solve the spectral and
  inverse spectral problem for the case of $m$ being a finite positive
  discrete measure. In particular, explicit determinantal formulas for
  the measure $m$ are given. These formulas generalize Stieltjes'
  formulas used by Krein in his study of the corresponding second
  order ODE $-\phi''=zm\phi$.
\end{abstract}



\submitto{\IP}


\maketitle

\section{Introduction}

In this paper we study the third order nonselfadjoint spectral problem
\begin{equation}
  \label{eq:cubic-neumann}
  \eqalign{
    -\phi'''(x)=zm(x)\phi(x) \quad\text{for $x\in\R$},\\
    \phi_x(-\infty)=\phi_{xx}(-\infty)=0,
    \qquad
    \phi_{xx}(\infty)=0,
  }
\end{equation}
which we refer to as the cubic string with Neumann-like
boundary conditions.
In particular we consider the case when $m(x)$ is a discrete measure
(a finite linear combination of Dirac deltas),
and the problem of reconstructing $m(x)$ given appropriate spectral data.
The purpose of this introduction is to explain where this problem comes from
and why it might be of interest.

Our starting point is the Degasperis--Procesi (DP) equation,
which is the completely integrable nonlinear PDE
\begin{equation}
  \label{eq:DP}
  u_t-u_{txx} + 4u u_x = 3 u_x u_{xx} + u u_{xxx}.
\end{equation}
It was discovered purely mathematically in a search for integrable equations
\cite{dp,dhh1},
but has been shown by Johnson \cite{johnson}
to play a role in shallow water wave theory similar to that of
the Camassa--Holm (CH) equation \cite{ch}
\begin{equation}
  \label{eq:CH}
  u_t-u_{txx} + 3u u_x = 2 u_x u_{xx} + u u_{xxx}.
\end{equation}
By taking appropriate scaling limits one can derive new nonlinear equations
which inherit some of the integrable structure of their ancestors.
Consider in particular the high-frequency limit
obtained by changing variables $x\mapsto \varepsilon x$, $t\mapsto \varepsilon t$,
and then letting $\varepsilon\to 0$.
For the CH equation \eref{eq:CH} this gives
the Hunter--Saxton equation \cite{hs,hz1} for nematic liquid crystals,
\begin{equation}
  \label{eq:HS}
  (u_t+u u_x)_{xx} = u_x u_{xx}.
\end{equation}
The DP equation \eref{eq:DP}, on the other hand, reduces in the high-frequency limit to
\begin{equation}
  \label{eq:dB}
  (u_t+u u_x)_{xx} = 0,
\end{equation}
which we will refer to as the derivative Burgers equation.
The name of course refers to the well known inviscid Burgers equation
(or Riemann shock wave equation)
$u_t+uu_x=0$,
the prototype equation for studying shock waves.
That \eref{eq:dB}, with the two extra derivatives, has a Lax pair
is a fairly recent discovery;
as far as we are aware, this was first shown by Qiao and Li \cite{qiao-li}
(who studied the whole integrable hierarchy that \eref{eq:dB} belongs to).

A unified way of writing these four equations is
\begin{equation}
  \label{eq:b-eqn}
  m_t + m_x u + b m u_x = 0,
\end{equation}
with $m$ and~$b$ as in the following table:
\begin{center}
  \begin{tabular}{c|c|c}
    & $m=u-u_{xx}$ & $m=u_{xx}$ \\
    \hline
    $b=2$ & Camassa--Holm & Hunter--Saxton \\
    \hline
    $b=3$ & Degasperis--Procesi & derivative Burgers\\
  \end{tabular}
\end{center}
(According to Dullin, Gottwald and Holm \cite{dullin},
the whole family \eref{eq:b-eqn} with $m=u-u_{xx}$ and $b\ne -1$
can be derived as an asymptotic limit of the shallow water equations.
The CH and DP equations are the only integrable cases in this family.
Moreover, \eref{eq:b-eqn} with $m=u-u_{xx}$
is included in the more general family
$u_t+[\delta u + 3\gamma u^2/2 + u^{1-\omega}(u^{\omega}u_x)_x ]_x + \nu u_{txx}=0$
studied by Rosenau \cite{rosenau};
take $1+2\omega=b$, $-3\gamma=b+1$, $\delta=0$, $\nu=-1$,
and let $x\mapsto -x$.
Another generalization of \eref{eq:b-eqn} is the pulson equation
introduced by Fringer and Holm \cite{fringer-holm-pulsons,holm-hone-pulsons}.)

A similar table summarizes the particular kinds of soliton-like
solutions characteristic of these equations, as well as the spectral
problems (including the type of boundary conditions)
employed when computing these solutions using inverse
spectral methods:
\begin{center}
  \begin{tabular}{c|c|c}
    & \textsc{Peakon solutions} & \textsc{Piecewise linear solutions} \\
    \hline
    \textit{Discrete string} & Camassa--Holm & Hunter--Saxton \\
    \hline
    \textit{Discrete cubic string} & Degasperis--Procesi & derivative Burgers\\
    \hline
    & \textit{Dirichlet(-like)} & \textit{Neumann(-like)} \\
  \end{tabular}
\end{center}
This requires perhaps some explanation.
Peakons (peaked solitons) are solutions taking the form
\begin{equation}
  u(x,t) = \sum_{k=1}^n m_k(t)\,\exp\bigl(- \abs{x-x_k(t)} \bigr),
\end{equation}
while the piecewise linear solutions are given by
\begin{equation}
  \label{eq:u-ansatz}
  u(x,t) = \sum_{k=1}^n m_k(t) \abs{x-x_k(t)},
\end{equation}
where in each case the functions $x_k(t)$ and $m_k(t)$ (positions and
momenta of the solitons)
are required to satisfy a certain system of $2n$ ODEs in order for
$u(x,t)$ to satisfy the PDE.
(See Theorem~\ref{thm:ODEs} below for the case of the
derivative Burgers equation.)

It was discovered in \cite{bss-accoustic} that the spatial equation
in the Lax pair of the CH equation can be transformed by a change
of variables to the string equation with Dirichlet boundary conditions,
\begin{equation}
  \label{eq:ordinarystring}
  \eqalign{
    -\phi''(x)=zm(x)\phi(x) \quad\text{for $-1<x<1$},\\
    \phi(-1)=\phi(1)=0,
  }
\end{equation}
where $z$ is the spectral parameter (squared frequency), $m$ is the
mass density and $\phi$ represents the amplitude of the harmonic mode
corresponding to to that frequency.
(This originates of course in separation of variables in the
linear wave equation $v_{xx}=m \, v_{tt}$ describing small vibrations
of a string with given mass density $m(x)$.)
The peakon solutions are related to the case when the mass distribution
is a sum of point masses (Dirac deltas);
this is what is meant by the \emph{discrete} string.
Krein \cite{krein-cord,krein-stieltjes,kackrein}
studied the inverse spectral problem of the string
equation \eref{eq:ordinarystring} and showed
in particular that in the discrete case the
mass distribution can be recovered from spectral data
using continued fractions of Stieltjes type \cite{stieltjes}.
This was exploited in \cite{bss-stieltjes,bss-moment}
to solve the CH peakon ODEs explicitly and analyze the solutions
in detail.
Similarly, the inverse problem for the string equation on $\R$
with Neumann boundary conditions
$\phi_x(-\infty)=\phi_x(\infty)=0$ was used in \cite{bss-hs}
to solve the ODEs for the piecewise linear solutions of the Hunter--Saxton
equation.
The inverse spectral theory for the simple looking ODE \eref{eq:ordinarystring} is
surprisingly rich in mathematical content. The most comprehensive
account of the inverse problem for this equation, and also of Krein's
work on the subject, is the monograph \cite{McKean} by Dym and McKean.

Despite the superficial similarity, the DP equation has a mathematical
structure quite different from that of the CH equation.
In particular,  the spatial equation in the Lax pair of the DP eqation is of third order
instead of second.
In \cite{ls-invprob} we showed how to transform it into the cubic string with
Dirichlet-like boundary conditions,
\begin{equation}
  \eqalign{
    -\phi'''(x)=zm(x)\phi(x) \quad\text{for $-1<x<1$},\\
    \phi(-1)=\phi_x(-1)=0,
    \qquad
    \phi(1)=0.
  }
\end{equation}
This problem was studied in great detail in our previous paper \cite{ls-imrp},
where the inverse problem was solved in the discrete case,
leading to explicit solution formulas for the DP peakons.

In this paper we first briefly show how the derivative Burgers
equation is related to the cubic string with Neumann-like
boundary conditions \eref{eq:cubic-neumann}.
Then we will give a fairly complete account of the direct and inverse
spectral theory of \eref{eq:cubic-neumann} in the case when
$m(x)$ is a positive discrete measure.
Some of the ideas used in \cite{ls-imrp} are developed further,
and we also provide new shorter proofs for some key steps which
perhaps appeared mysterious there. For the most part this
paper is self-contained, with the exception of using some Heine-like
multiple integral formulas proved in \cite{ls-imrp} (see Appendix~B)
and postponing the proof of Theorem~\ref{thm:4recurrence} to another paper.

\section{The derivative Burgers equation}

Here we give a short account of the basic facts about the derivative Burgers equation
\eref{eq:dB}.
First we show how it inherits the integrable structure of the DP
equation \eref{eq:DP}. The following Lax pair for the DP equation was
given in \cite{dhh1}:
\begin{equation}
  \label{eq:DPlax}
  (\partial_x - \partial_x^3) \phi = z m\phi,
  \qquad
  \phi_t = \left[ z^{-1} (c-\partial_x^2) + u_x - u \partial_x \right] \phi,
\end{equation}
where $c$ is an arbitrary constant (or function of~$z$);
indeed, the compatibility condition $\phi_{txxx}=\phi_{xxxt}$
is equivalent to $m_t+m_x u+3m u_x=0$ and $m_x=(u-u_{xx})_x$.
The high-frequency substitution $x\mapsto \varepsilon x$, $t\mapsto \varepsilon t$,
together with $m\mapsto -\varepsilon^{-2}m$, $z\mapsto -\varepsilon^{-1}z$,
and $c\mapsto -z \varepsilon^{-2}c$,
yields
\begin{equation*}
  (\varepsilon^2 \partial_x - \partial_x^3) \phi = z m\phi,
  \qquad
  \phi_t = \left[ z^{-1} \partial_x^2 + c + u_x - u \partial_x \right] \phi,
\end{equation*}
which is compatible iff $m_t+m_x u+3m u_x=0$ and $m_x=(u_{xx}-\varepsilon^2 u)_x$.
Letting $\varepsilon\to 0$ we obtain the system
\begin{eqnarray}
  \label{eq:lax1}
  - \partial_x^3 \phi = z m\phi, \\
  \label{eq:lax2}
  \phi_t = \left[ z^{-1} \partial_x^2 + c + u_x - u \partial_x
  \right] \phi,
\end{eqnarray}
which is compatible iff $m_t+m_x u+3m u_x=0$ and $m_x=u_{xxx}$,
and hence in particular when the derivative Burgers equation holds.

\begin{remark}
  Here the cubic string \eref{eq:lax1} on the whole
  real line $\R$ appears naturally. This is in contrast to the DP case
  \cite{ls-invprob} where a change of variables mapping the real line
  to the finite interval $(-1,1)$ is required to get rid of the
  $\partial_x$ term in the first equation of \eref{eq:DPlax}.
  A similar distinction holds for the ordinary string equation
  in the Hunter--Saxton and CH cases.
\end{remark}

\begin{remark}
  Hone and Wang \cite{hone-wang} have described a different
  high-frequency limititing procedure taking the DP equation to
  (the $x$ derivative of) the integrable Vakhnenko equation
  $(u_t+uu_x)_x+u=0$ \cite{vakhnenko,vakhnenko-parkes-DP-V-eqn}.
  It should perhaps also be mentioned that in the low-frequency limit
  $\epsilon\to\infty$ the substitution above directly reduces the
  DP equation to $u_t+uu_x=0$ \cite{dhh1},
  but this does not give a useful Lax pair.
\end{remark}

We will return to the question of boundary conditions for
\eref{eq:lax1} shortly.
The piecewise linear solutions of the derivative Burgers equation
are described by the following theorem,
where we use the convention $\sgn 0 = 0$.
Dots denote $\frac{\rmd}{\rmd{}t}$ as usual.

\begin{theorem}
  \label{thm:ODEs}
  The function $u$ given by \eref{eq:u-ansatz} satisfies the derivative Burgers
  equation \eref{eq:dB} in the sense of distributions
  if and only if
  \begin{equation}
    \label{eq:ODEs}
     \dot{x}_k =\sum_{i=1}^nm_i\abs{ x_k-x_i},
    \qquad
    \dot{m}_k = 2 \sum_{i=1}^n m_k m_i \sgn(x_i-x_k)
  \end{equation}
  for $k=1,\ldots,n$.
\end{theorem}

One can assume that all $m_k\neq 0$,
since it follows from \eref{eq:ODEs} that any vanishing $m_k$ remains zero
and hence never enters the solution.

The solutions in Theorem~\ref{thm:ODEs} satisfy the usual Burgers
equation $u_t+uu_x=0$ (without the extra derivatives) only if $\sum_{k=1}^n m_k=0$.
This case will not be considered here, since we will make the
following basic assumptions throughout this paper:
\begin{equation}
  \label{eq:ordering-positive}
  x_1<\cdots<x_n,\qquad \text{$m_k>0$ for all $k$}.
\end{equation}
It can be shown that if these two assumptions hold for $t=0$,
then they will continue to hold, and the solution of \eref{eq:ODEs} will be defined,
for all $t>0$.
(But if $m_k$'s of both signs are present,
shocks may form after finite time, and then the solution will no longer
have the form \eref{eq:u-ansatz}.)
Under the assumptions \eref{eq:ordering-positive},
$\sgn(x_i-x_k)$ can be replaced by $\sgn(i-k)$ in \eref{eq:ODEs}.

\begin{proposition}
  \label{prop:M-Mplus}
  $M=\sum_{k=1}^n m_k$ and
  $M_+=\sum_{k=1}^n m_k x_k$ are constants of motion of
  the ODEs \eref{eq:ODEs}.
\end{proposition}

\begin{proof}
  This is easily verified by a short calculation.
\end{proof}

Since we are now in the context of piecewise linear solutions 
where $u$ has the form \eref{eq:u-ansatz},
it follows that $m=u_{xx}=2\sum_1^n m_k\,\delta(x-x_k)$
is a (time-dependent) discrete measure.
Consider some fixed time~$t$.
The wave function $\phi(x,t;z)$ satisfying the cubic string equation
\eref{eq:lax1} is piecewise a quadratic polynomial in~$x$,
since $\phi_{xxx}=0$ away from the support of~$m$.
At the points $x_k$, \eref{eq:lax1} implies that $\phi$ and $\phi_x$ are continuous
while $\phi_{xx}$ has a jump of size $-2 z\,m_k \,\phi(x_k)$.
These requirements completely define $\phi$
once $\phi$, $\phi_x$, and $\phi_{xx}$ are prescribed at some point.

Now let $u$ evolve according to the derivative Burgers equation;
that is, let $x_k$ and $m_k$ evolve according to \eref{eq:ODEs}.
Then $\phi$ evolves according to \eref{eq:lax2}.
We can postulate that
\begin{equation}
  \label{eq:phileft}
  \phi(x,t;z) = 1 \qquad \text{for} \quad x < x_1(t)
\end{equation}
provided that we choose $c$ in \eref{eq:lax2}
equal to the constant of motion~$M=\sum m_k$;
then \eref{eq:phileft} is compatible with the time evolution
since it makes both sides of \eref{eq:lax2} vanish for $x<x_1$.
Propagating \eref{eq:phileft} to the right gives
\begin{equation}
  \label{eq:phiright}
  \phi(x,t;z) = A(t;z)\,\frac{(x-x_n)^2}{2} + B(t;z)\,(x-x_n) + C(t;z)
\end{equation}
for $x > x_n(t)$,
where $A$, $B$, $C$ are some explicitly computable polynomials in~$z$ of
degree~$n$,
with coefficients depending on the $x_k(t)$'s and $m_k(t)$'s;
see \eref{eq:ABC} below.
Inserting this into \eref{eq:lax2}, together with
$u=\sum m_k(x-x_k)=Mx-M_+$ for $x>x_n(t)$,
gives the time evolution of these quantities:
\begin{equation}
  \label{eq:ABCdot}
    \dot{A} = 0,
    \qquad
    \dot{B} =MB,
    \qquad
    \dot{C} = \frac{A}{z} + 2MC.
\end{equation}
We see that for $x>x_n$ we can impose $\phi_{xx}=2A=0$; this
condition will be preserved by the evolution in time. Together
with \eref{eq:phileft}, this amounts to $\phi$ satisfying
the Neumann-like boundary conditions
$\phi_x(-\infty)=\phi_{xx}(-\infty)=0$, $\phi_{xx}(\infty)=0$,
which is how we were naturally led to the cubic string
boundary value problem \eref{eq:cubic-neumann}.

More details about the solutions of the ODEs \eref{eq:ODEs}
will be given in a separate paper
devoted to the derivative Burgers equation.
Here we will leave this subject after showing that
$M$ is the first in a naturally appearing sequence of
constants of motion,
while $M_+$ is something of an odd bird.
Indeed, writing $A(z) = 2 \sum_{k=1}^n (-z)^k M_k$
it is straightforward to derive an expression for the coefficients
$M_k$ using \eref{eq:ABC},
and since $A$ is constant in time, we then obtain the following result.

\begin{theorem}
  \label{thm:constants-of-motion}
  The functions $M_1,\ldots,M_n$ given by
  \begin{equation}
    \label{eq:Mk}
    M_k =
    \sum_{I\in\binom{[1,n]}{k}}
    \Bigl(
    \prod_{i \in I} m_i
    \Bigr)
    \Bigl(
    \prod_{j=1}^{k-1} (x_{i_j}-x_{i_{j+1}})^2
    \Bigr).
  \end{equation}
  are constants of motion for the ODEs \eref{eq:ODEs}.
  Here $\binom{[1,n]}{k}$ denotes the set of all $k$-element subsets
  $I=\left\{ i_1 < \cdots < i_k \right\}$ of the integer interval
  $[1,n]=\left\{ 1,\ldots,n \right\}$.
  The empty product when $k=1$ is to be interpreted as $1$, so that $M_1=\sum m_k=M$.
\end{theorem}

For example, when $n=3$ the constants of motion are
\begin{equation}
  \eqalign{
    M_1 &= m_1+m_2+m_3 = M,
    \\
    M_2 &= m_1 m_2 (x_1-x_2)^2 + m_1 m_3 (x_1-x_3)^2 + m_2 m_3 (x_2-x_3)^2,
    \\
    M_3 &= m_1 m_2 m_3 (x_1-x_2)^2 (x_2-x_3)^2,
  }
\end{equation}
together with $M_+ = m_1 x_1+m_2 x_2+m_3 x_3$.

\section{The forward spectral problem}
\label{sec:spectral}

We now turn to the spectral theory of the
cubic string with Neumann-like boundary conditions
\eref{eq:cubic-neumann}.
The problem makes sense for any $m(x)$ decaying sufficiently fast as $\abs{x}\to\infty$,
but we will deal exclusively with the discrete case
where $m(x)=2\sum m_k \, \delta(x-x_k)$,
and in addition we assume ordering and positivity as in
\eref{eq:ordering-positive}.
Then the eigenvalues $z$ are the zeros of the time-independent
polynomial~$A(z)$ defined by \eref{eq:phiright},
so the ODEs \eref{eq:ODEs} associated to the
derivative Burgers equation induce an isospectral deformation of
\eref{eq:cubic-neumann}.
If we suppress the $t$ dependence in the notation,
the propagation of the wave function from left to right is described explicitly by
\begin{equation}
  \label{eq:propagation}
  \Phi(x_{k+1}^-) = L_k \, \Phi(x_k^+),
  \qquad
  \Phi(x_k^+) = G_k(z) \, \Phi(x_k^-),
\end{equation}
where $\Phi=(\phi,\phi_x,\phi_{xx})^t$ and,
with $l_k=x_{k+1}-x_k$,
\begin{equation}
  L_k =
  \begin{pmatrix}
    1 & l_k & l_k^2 / 2 \\
    0 & 1 & l_k \\
    0 & 0 & 1
  \end{pmatrix},
  \qquad
  G_k(z) =
  \begin{pmatrix}
    1 & 0 & 0 \\
    0 & 1 & 0 \\
    -2 m_k z & 0 & 1
  \end{pmatrix}.
\end{equation}
Hence
\begin{equation}
  \label{eq:ABC}
  \fl
  \begin{pmatrix}
    C(z) \\ B(z) \\ A(z)
  \end{pmatrix}
  = \Phi(x_n^+;z)
  = G_n(z) L_{n-1} G_{n-1}(z) \ldots L_2 G_2(z) L_1 G_1(z)
  \begin{pmatrix}
    1 \\ 0 \\ 0
  \end{pmatrix}.
\end{equation}
From this one can easily extract the coefficients of $A(z) = 2 \sum_{k=1}^n (-z)^k M_k$
explicitly; see Theorem~\ref{thm:constants-of-motion} above.
Note also that, as $z\to 0$,
\begin{equation}
  \label{eq:BCfirstterms}
  \eqalign{
    B(z) &= 2z (M_+ - Mx_n) + \Or(z^2),
    \\
    C(z) &= 1- z \sum_{k=1}^n m_k (x_n-x_k)^2 + \Or(z^2).
  }
\end{equation}

\begin{theorem}
  \label{thm:real-spectrum}
  If $m_k>0$ for all~$k$, then the spectral problem \eref{eq:cubic-neumann}
  with $m(x)=2\sum m_k \, \delta(x-x_k)$
  has $n$ distinct nonnegative eigenvalues $z=\lambda_k$,
  \begin{equation}
    0=\lambda_0 < \lambda_1 < \cdots < \lambda_{n-1}.
  \end{equation}
\end{theorem}

\begin{proof}
  It is trivial that $z=0$ is an eigenvalue, with corresponding
  eigenfunction $\phi=1$. The proof of the other assertions amounts to
  showing that the remaining eigenvalues coincide with the eigenvalues
  of a certain oscillatory matrix. The basic facts about
  oscillatory matrices are summarized in Appendix~A.
  By \eref{eq:propagation}, we have $\Phi(x_{k+1}^+)=G_{k+1}(z) \, L_k \, \Phi(x_k^+)$.
  With $(C_k,B_k,A_k)^t=\Phi(x_k^+)$, this says that
  \begin{equation*}
    \eqalign{
      C_{k+1} &= C_k + l_k B_k + \frac{l_k^2}{2} A_k,\\
      B_{k+1} &= B_k + l_k A_k,\\
      A_{k+1} &= -2z\,m_{k+1} C_{k+1} + A_k.
    }
  \end{equation*}
  The boundary conditions at $-\infty$ can be expressed as $A_0=B_0=0$.
  The first equation gives
  \begin{equation*}
    2z\,C_{k+1} -2z\,C_k = 2z\,l_k B_k + z\,l_k^2 A_k,
  \end{equation*}
  which, using the second and third equations, becomes
  \begin{equation*}
    -\frac{A_{k+1}-A_k}{m_{k+1}} + \frac{A_k-A_{k-1}}{m_k}
    = 2z\,l_k \sum_{j=1}^{k-1} l_j A_j + z\,l_k^2 A_k.
  \end{equation*}
  The boundary condition at $+\infty$ is $A=A_n=0$,
  so $\matrixX=(A_1,\ldots,A_{n-1})^t$ satisfies
  $\matrixM \matrixX=z\matrixL \matrixX$,
  where $\matrixM$ is the symmetric tridiagonal matrix defined by
  $\matrixM_{ii}=m_i^{-1}+m_{i+1}^{-1}$
  and $\matrixM_{i,i-1}=\matrixM_{i-1,i}=-m_i^{-1}$,
  while $\matrixL$ is the lower triangular matrix with
  $\matrixL_{ii}=l_i^2$ and $\matrixL_{ij}=2 l_i l_j$ for $i>j$.
  ($\matrixM$ and $\matrixL$ are both of size $(n-1)\times(n-1)$.)
  An easy induction on~$n$ shows that
  \begin{equation*}
    \det \matrixM = \frac{m_1+m_2+\cdots+m_n}{m_1 m_2 \ldots m_n},
  \end{equation*}
  which is positive when all $m_k>0$, so $\matrixM^{-1}$ exists.
  Hence the nonzero eigenvalues of \eref{eq:cubic-neumann}
  are the reciprocals of the eigenvalues of the matrix $\matrixM^{-1}\matrixL$.

  By Theorem~\ref{thm:TPspectrum} we are done if we can show that
  $\matrixM^{-1}\matrixL$ is oscillatory.
  The star operation does not alter the determinant,
  since $\matrixM^*=\matrixT \matrixM \matrixT$ with
  $\matrixT=\diag(+1,-1,+1,\ldots,\pm1)$.
  Hence the determinant of $\matrixM^*$ is positive, and so are also all its
  principal minors, since they are determinants of the same form.
  Thus $\matrixM^*$ is positive definite,
  and Theorem~\ref{thm:tridiagonal} then implies that $\matrixM^*$ is oscillatory.
  Hence $\matrixM^{-1}=((\matrixM^*)^*)^{-1}$ is oscillatory by Theorem~\ref{thm:signoscillatory}.

  \begin{figure}
    \centering
    \begin{equation*}
      \xymatrix{
        \text{\small Source $1$} & \bullet \ar[r]^{l_1} & \bullet \ar[r] & \bullet \ar[r]^{l_1} & \bullet & \text{\small Sink $1$} \\
        \text{\small Source $2$} & \bullet \ar[r]^{l_2} & \bullet \ar[r] \ar[u] \ar[ur] & \bullet \ar[r]^{l_2} & \bullet & \text{\small Sink $2$} \\
        & \bullet \ar[r]^{l_3} & \bullet \ar[r] \ar[u] \ar[ur] & \bullet \ar[r]^{l_3} & \bullet & \\
        \raisebox{0pt}[0pt][0pt]{\small\vdots} & \ar@{.>}[r] &  \ar@{.>}[r] \ar@{.>}[u] \ar@{.>}[ur] &  \ar@{.>}[r] & & \raisebox{0pt}[0pt][0pt]{\small\vdots} \\
        \text{\small Source $n-1$} & \bullet \ar[r]^{l_{n-1}} & \bullet \ar[r] \ar[u] \ar[ur] & \bullet \ar[r]^{l_{n-1}} & \bullet & \text{\small Sink $n-1$}
      }
    \end{equation*}
    \caption{A planar network of order $n-1$ whose weighted path matrix is $\matrixL$. Unlabelled edges have weight~$1$.}
    \label{fig:network}
  \end{figure}
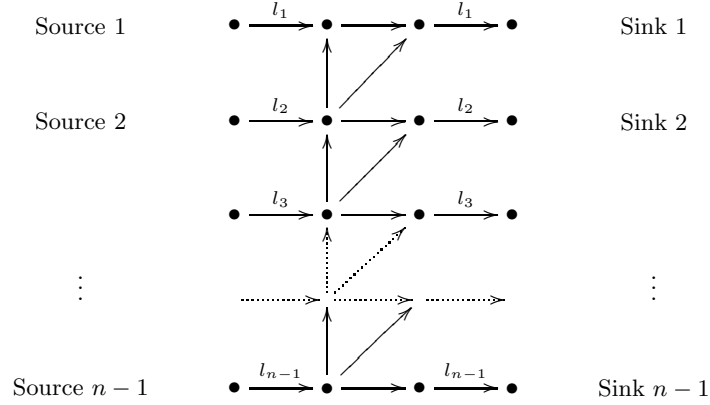

  The matrix $\matrixL$ is the weighted path matrix of the planar network in
  Figure~\ref{fig:network},
  and consequently totally nonnegative by Theorem~\ref{lem:Lind}
  (Lindstr\"om's Lemma).
  Since $\matrixL$ is also obviously nonsingular,
  $\matrixM^{-1}\matrixL$ is oscillatory by Theorem~\ref{thm:product}.
  The proof is finished.
\end{proof}

Next, we define the \emph{Weyl functions} $W(z)$ and $Z(z)$ of the
discrete cubic string with Neumann-like boundary conditions, as well
as the residues $b_k$ and $c_k$ in their partial fraction
decompositions.
The case $M=\sum m_k=0$
which would need special treatment is excluded here by our
standing assumption that all $m_k>0$.
Recall that the zeros of $A(z)$ are $z=0$ and $z=\lambda_k$,
and that the lowest order terms of $B(z)$ and $C(z)$ are given
by \eref{eq:BCfirstterms}.

\begin{definition}
  \label{def:WZdef}
  Let
  \begin{equation}
    \label{eq:WZdef}
    \eqalign{
      W(z) = \frac{\phi_x(x_n+)}{\phi_{xx}(x_n+)}
      = \frac{B(z)}{A(z)}
      = \sum_{k=1}^{n-1} \frac{b_k}{z-\lambda_k},
      \\
      Z(z) = \frac{\phi(x_n+)}{\phi_{xx}(x_n+)}
      = \frac{C(z)}{A(z)}
      = -\frac{1/2M}{z} + \sum_{k=1}^{n-1} \frac{c_k}{z-\lambda_k}.
    }
  \end{equation}
\end{definition}

\begin{theorem}
  \label{thm:WdeterminesZ}
  The Weyl functions satisfy the relation
  \begin{equation}
    \label{eq:ZW-relation}
    Z(z) + Z(-z) = W(z)W(-z).
  \end{equation}
  Hence, the second Weyl function $Z$ is determined by the first Weyl function $W$
  (except for the residue at $z=0$ which is determined by the constant $M$)
  through the formula
  \begin{equation}
    \label{eq:cb-relation}
    c_k = W(-\lambda_k) \, b_k = -\sum_{j=1}^{n-1} \frac{b_j b_k}{\lambda_j + \lambda_k},
    \qquad k=1,\ldots,n-1.
  \end{equation}
\end{theorem}

\begin{proof}
  Let
  \begin{equation}
    \label{eq:transition-matrix}
    S(z) = G_n(z) L_{n-1} G_{n-1}(z) \ldots L_2 G_2(z) L_1 G_1(z)
  \end{equation}
  denote the transition matrix from \eref{eq:ABC},
  whose first column is $(C,B,A)^t$.
  The $(1,3)$ entry of the identity
  \begin{equation}
    \label{eq:automorphism}
    S(-z)^t J S(z) J = I,
    \qquad\text{where}\quad
    J =
    \begin{pmatrix}
      0 & 0 & 1 \\
      0 & -1 & 0 \\
      1 & 0 & 0
    \end{pmatrix}
  \end{equation}
  reads
  \begin{equation*}
    A(-z) C(z) - B(-z)B(z) + C(-z) A(z) = 0.
  \end{equation*}
  (Equation \eref{eq:automorphism} is proved by noting that $X=G_k(z)$
  and $X=L_k$ both satisfy $X(-z)^{-1}=(JX(z)J)^t$,
  therefore so does $X=S(z)$ \cite{ls-imrp}.)
  By definition, $W=B/A$ and $Z=C/A$,
  so division by $A(z)A(-z)$ yields \eref{eq:ZW-relation}.
  The residue of \eref{eq:ZW-relation} at $z=\lambda_k$ is \eref{eq:cb-relation}.
\end{proof}

\begin{remark}
  This proof (with trivial modifications) also works for proving the
  corresponding relation between the two Weyl functions for the
  the discrete cubic string with Dirichlet-like boundary conditions,
  which is used in the DP case.
  It is more straightforward than the proof in \cite{ls-imrp},
  which relied on the analogue of Theorem~\ref{thm:weyl-evolution} below.
\end{remark}

\begin{theorem}
  \label{thm:residues}
  If $m_k>0$ for all~$k$, then the residues $b_k$ of the Weyl function $W(z)$
  are all negative.
\end{theorem}

\begin{proof}
  By definition $A$, $B$, and $C$ are all positive when $z<0$.
  Hence $W(-\lambda_k)>0$, so that $b_k=B(\lambda_k)/A'(\lambda_k)$
  and $c_k=C(\lambda_k)/A'(\lambda_k)$ are of the same sign
  by \eref{eq:cb-relation};
  they cannot be zero, since that would imply that $A=B=C=0$ for $z=\lambda_k$,
  contradicting $\det S(z)=1$.

  We now use induction on $n$.
  For $n=2$ the Weyl function is
  \begin{equation*}
    \fl
    W(z)
    = \frac{B(z)}{A(z)}
    = \frac{-2 m_1 l_1 \,z}{-2(m_1+m_2)z+2(m_1 m_2 l_1^2)z^2}
    = \frac{-1/m_2 l_1}{z - \dsty\frac{m_1+m_2}{m_1 m_2 l_1^2}}
    = \frac{b_1}{z-\lambda_1},
  \end{equation*}
  so $b_1<0$.
  For $n=p>1$, if we let $m_1\to 0^+$ then $W=W_p$ tends to the
  corresponding Weyl function $W_{p-1}$ for the cubic string with parameters
  $m_2,\ldots,m_p$ and $l_2,\ldots,l_{p-1}$.
  For the largest eigenvalue and corresponding residue,
  we have $\lambda_{p-1} \to +\infty$ and $b_{p-1} \to 0$,
  while the other eigenvalues and residues tend to the eigenvalues and
  (negative, by induction hypothesis) residues of $W_{p-1}$.
  Since the $\lambda_k$'s and $b_k$'s depend continuously on the parameters,
  the $b_k$'s cannot change sign during this process,
  and thus they must have been negative to begin with,
  except perhaps for $b_{p-1}$.

  To prove that $b_{p-1}<0$, note that $B$, by \eref{eq:ABC} and
  the definition of $G_p(z)$ and $L_{p-1}$, does not depend on $m_p$.
  By making $m_p$ small, we can make $\lambda_{p-1}$ as large as desired;
  in particular, larger in absolute value than all zeros of $B$,
  which except for $z=0$ are also the zeros of $W$.

  Since the highest coefficients of $A$ and $B$ differ in sign,
  $W<0$ for all sufficiently large $z>0$.
  If $b_{p-1}$ were positive, then $W$ would be positive for some $z$ slightly
  larger than~$\lambda_{p-1}$, and that would force $W$ to have a zero
  which is larger than~$\lambda_{p-1}$, in contradiction to what we found above.
  Thus $b_{p-1}<0$.
  \end{proof}

\begin{remark}
  Rather than using one spectrum of eigenvalues together with
  the corresponding residues,
  it is clear from \eref{eq:WZdef}
  that the Weyl function $W(z)$ could instead be defined
  by two spectra, namely the zeros of $\phi_x(x_n+)$ and $\phi_{xx}(x_n+)$.
  Theorem~\ref{thm:residues} above implies that these two sets of zeros
  are interlacing.
  In this regard, and also in
  being determined by a single Weyl function (which is perhaps surprising),
  the cubic string is similar to the classical Sturm--Liouville problem and
  to the ordinary string.
\end{remark}

Although it is not needed for the purposes of this paper,
we also note how the Weyl functions change when the configuration
$\{ x_k,m_k \}_{k=1}^n$ of the discrete cubic string
evolves in time according to the ODEs \eref{eq:ODEs}.

\begin{theorem}
  \label{thm:weyl-evolution}
  The time evolution of the Weyl functions is given by
  \begin{equation}
    \dot{W}(z) = M \, W(z),
    \qquad
    \dot{Z}(z) = \frac{1}{z} + 2M \, Z(z).
  \end{equation}
  Hence, the residues evolve according to
  \begin{equation}
    b_k(t) = b_k(0) \, \rme^{Mt},
    \qquad
    c_k(t) = c_k(0) \, \rme^{2Mt}.
  \end{equation}
\end{theorem}

\begin{proof}
  This follows immediately from \eref{eq:ABCdot}.
\end{proof}

\section{The inverse spectral problem}
\label{sec:inverse-spectral}

A given discrete cubic string configuration $\{x_k,m_k\}_{k=1}^n$
with $x_1<\cdots<x_n$ and all $m_k>0$ defines, as described in the previous
section, a set of spectral data consisting of the positive eigenvalues
$\lambda_1<\cdots<\lambda_{n-1}$ together with their corresponding residues
$b_1,\ldots,b_{n-1}<0$ and the constant of motion $M=\sum m_k>0$.
Only the relative positions $l_k=x_{k+1}-x_k$ enter this computation,
which is why $2n-1$ is the correct number of spectral data if we wish
to reconstruct the $2n-1$ string parameters $l_1,\ldots,l_{n-1}$ and
$m_1,\ldots,m_n$.

Let $a^{(m)}$ denote the product of the first $m$ factors in the transition matrix $S(z)$
of \eref{eq:transition-matrix}:
\begin{equation}
  \label{eq:k-step-matrix}
  \fl
  \eqalign{
    a^{(2k)}(z)
    = G_n(z) L_{n-1} G_{n-1}(z) L_{n-2} \ldots G_{n-k+1}(z) L_{n-k},
    \\
    a^{(2k+1)}(z)
    = G_n(z) L_{n-1} G_{n-1}(z) L_{n-2} \ldots G_{n-k+1}(z) L_{n-k} G_{n-k}(z).
}
\end{equation}
The Weyl functions $W=B/A$ and $Z=C/A$ are constructed from the
first column $\Phi(x_n^+)=(C,B,A)^t$ in $S(z)=a^{(2n-1)}$.
From the entries in the odd partial products $a^{(2k+1)}(z)$, which are polynomials in~$z$,
we will show how to construct rational approximations to $W$ and $Z$
with a common denominator, and satisfying certain normalization and symmetry properties
rendering the approximants unique.
Conversely, given a set of spectral data, we define $W$ and $Z$ by \eref{eq:WZdef}.
From the set of uniquely determined rational approximants referred to above
we reconstruct all the factors $G_j$ and $L_j$ and
thus solve the inverse problem.

First, we have a simple power counting lemma.
Recall that $\Phi$ denotes the vector $(\phi,\phi_x,\phi_{xx})$.
For $0\le k \le n-1$, let
\begin{equation}
  \label{eq:pqr}
  \begin{pmatrix}
    p_k(z) \\ q_k(z) \\ r_k(z)
  \end{pmatrix}
  = \Phi(x_{n-k}^-;z),
\end{equation}
so that $(C,B,A)^t = \Phi(x_n^+) = a^{(2k+1)} (p_k,q_k,r_k)^t$
for any~$k$.

\begin{lemma}
  \label{lem:a-deg}
  If we compute the polynomial degree of a matrix entry-wise, then
  \begin{equation}
    \degree a^{(2k+1)}(z)=\begin{pmatrix}k&k-1&k-1\\k&k-1&k-1\\k+1&k&k \end{pmatrix}
  \end{equation}
  with the proviso that for $k=0$ the degrees in the second and third
  columns are all zero (negative degrees counted as zero).  Furthermore,
  $\degree(p_k)=\degree(q_k)=\degree(r_k)=n-k-1$.
\end{lemma}

\begin{proof}
  Since each $G_j(z)$ is linear in~$z$, and since the leftmost factor $G_n(z)$
  in \eref{eq:k-step-matrix}
  only affects the third row of $a^{(2k)}(z)$,
  the first two rows of $a^{(2k)}(z)$ have degree \mbox{$k-1$} while
  the third row has degree~$k$.
  Subsequent multiplication by $G_{n-k}(z)$ on the right gives the first claim.
  The second claim follows from an easy induction
  on $k'=n-k$ which starts with
  $(p_{n-1},q_{n-1},r_{n-1})=(1,0,0)$ for~$k'=1$.
\end{proof}

Fix $k$ and write $a(z)$ instead of $a^{(2k+1)}(z)$ for simplicity.
We now turn our attention to approximation properties of the entries of
$a(z)$.  These properties will later model the approximation problems
required to solve the inverse problem.

\begin{theorem}
  \label{thm:a-approx}
  For $k$ fixed, and $z\to\infty$,
  \begin{equation}
    \label{eq:WZapprox}
    \eqalign{
      a_{31} \begin{pmatrix} Z\\W \end{pmatrix}
      - \begin{pmatrix} a_{11}\\a_{21} \end{pmatrix}
      & =\begin{pmatrix} \Or(z^{-1})\\\Or(1) \end{pmatrix},
      \\
      a_{32}\begin{pmatrix} Z\\W \end{pmatrix}
      - \begin{pmatrix} a_{12}\\a_{22} \end{pmatrix}
      & =\begin{pmatrix} \Or(z^{-1})\\\Or(1) \end{pmatrix},
      \\
      a_{33} \begin{pmatrix} Z\\W \end{pmatrix}
      - \begin{pmatrix} a_{13}\\a_{23} \end{pmatrix}
      & =\begin{pmatrix} \Or(z^{-2})\\\Or(z^{-1}) \end{pmatrix},
    }
  \end{equation}
  where
  \begin{equation}
    \label{eq:WZnormalization}
    \eqalign{
      \begin{pmatrix} a_{11}(0) \\ a_{21}(0) \end{pmatrix}
      = \begin{pmatrix} 1\\0 \end{pmatrix},
      \qquad
      \begin{pmatrix} a_{12}(0) \\ a_{22}(0) \end{pmatrix}
      = \begin{pmatrix} \sum_{j=n-k}^{n-1} l_j \\1 \end{pmatrix},
      \\
      \ms
      a_{31}(0)=a_{32}(0)=0, \qquad a_{33}(0)=1.
    }
  \end{equation}
  Moreover, setting $W^*(z)=-W(-z)$ and $Z^*(z)=Z(-z)$ we have
  \begin{equation}
    \label{eq:WZsymmetry}
    a_{3j} Z^*+ a_{2j} W^* +a_{1j} = \Or(z^{-(k+1)}),
    \qquad 1\le j\le 3.
  \end{equation}
\end{theorem}

\begin{proof}
  Regarding \eref{eq:WZapprox}
  we only give the proof for the third column of $a(z)$ to illustrate the main steps.
  Write $(p,q,r)$ for $(p_k,q_k,r_k)$.
  Since
  \begin{equation*}
    \begin{pmatrix}
      Z \\ W \\ 1
    \end {pmatrix}
    = \frac{1}{A}
    \begin{pmatrix}
      C \\ B \\ A
    \end{pmatrix}
    = \frac{1}{A} \, a(z)
    \begin{pmatrix}
      p \\ q \\ r
    \end{pmatrix},
  \end{equation*}
  we have
  \begin{equation*}
    \fl
    W - \frac{a_{23}}{a_{33}}
    = \frac{a_{21}p+a_{22}q+a_{23}r}{a_{31}p+a_{32}q+a_{33}r} - \frac{a_{23}}{a_{33}}
    = \frac{\dsty -\frac{a_{31}a_{23}-a_{21}a_{33}}{a_{33}} \frac{p}{r}
      +\frac{a_{22}a_{33}-a_{32}a_{23}}{a_{33}} \frac{q}{r}}
    {a_{31}\frac{p}{r}+a_{32}\frac{q}{r}+a_{33}}.
  \end{equation*}
  The $2\times2$ minors of $a(z)$ appearing here are the $(2,1)$ and $(1,1)$
  entries of the adjoint matrix $\adj\bigl(a(z)\bigr)$, which equals $a(z)^{-1}$ since
  each factor in $a(z)$ has determinant~$1$.
  Hence the matrix
  \begin{equation*}
    \fl
    \adj\bigl(a(z)\bigr) =
    G_{n-k}(-z) L_{n-k}^{-1} G_{n-k+1}(-z) \ldots L_{n-2}^{-1} G_{n-1}(-z) L_{n-1}^{-1} G_n(-z)
  \end{equation*}
  has the same degree structure as in Lemma~\ref{lem:a-deg},
  so the minors above both have degree~$k$, and so do $a_{32}$ and $a_{33}$.
  Moreover, $p/r$ and $q/r$ are both of order $\Or(1)$ by Lemma~\ref{lem:a-deg}, while
  $a_{31}$ has degree \mbox{$k+1$},
  so the whole expression is of order $\Or(z^{-(k+1)})$ as claimed.
  The proof for $Z-a_{13}/a_{33}$ is similar, except that the appearing minors
  come from the first two rows of the second column of $\adj(a(z))$ instead,
  which lowers the degree one step.

  \Eref{eq:WZnormalization} follows from \eref{eq:k-step-matrix}
  since $G_j(0)=I$.
  As for the symmetry relation \eref{eq:WZsymmetry},
  \begin{equation*}
    a(z)^{-1}
    \begin{pmatrix}
      Z \\ W \\ 1
    \end{pmatrix}
    = \frac{1}{A}
    \begin{pmatrix}
      p \\ q \\ r
    \end{pmatrix}
    =
    \begin{pmatrix}
      \Or(z^{-(k+1)}) \\
      \Or(z^{-(k+1)}) \\
      \Or(z^{-(k+1)})
    \end{pmatrix}.
  \end{equation*}
  Like $S(z)$, the matrix $a(z)$ satisfies the property \eref{eq:automorphism}
  (the proof is the same),
  so the left-hand side above equals
  $J a(-z)^t J \, (Z,W,1)^t$,
  the rows of which give rise to \eref{eq:WZsymmetry}
  after changing $z$ to~$-z$.
\end{proof}

\begin{proposition}
  For $k=0$,
  \begin{equation}
    \label{eq:kzerocase}
    a^{(1)}_{31}(z)=-2m_nz, \quad a^{(1)}_{32}(z)=0, \quad a^{(1)}_{33}(z)=1.
  \end{equation}
  For $k=1,\ldots,n-1$,
  \begin{eqnarray}
    \fl
    a^{(2k+1)}_{31}(z) = -2z\left( \sum_{i=n-k}^{n} m_i \right)
    + \cdots +  (-2z)^{k+1} \left( \prod_{i=n+1-k}^{n} m_i \frac{l_{i-1}^2}{2} \right)
    m_{n-k},
    \label{eq:a31}
    \\
    \fl
    a^{(2k+1)}_{32}(z) = -2z \left( \sum_{i=n+1-k}^{n} m_i (x_i-x_{n-k}) \right)
    \nonumber \\
    + \cdots
    + (-2z)^{k} \left( \prod_{i=n+2-k}^{n} m_i \frac{l_{i-1}^2}{2} \right) m_{n+1-k} l_{n-k},
    \label{eq:a32}
    \\
    \fl
    a^{(2k+1)}_{33}(z) = 1 - 2z \left( \sum_{i=n+1-k}^n m_i \frac{(x_i-x_{n-k})^2}{2} \right)
    \nonumber \\
    + \cdots
    + (-2z)^k \left( \prod_{i=n+1-k}^n m_i \frac{l_{i-1}^2}{2} \right).
    \label{eq:a33}
  \end{eqnarray}
\end{proposition}

\begin{proof}
  Write $G_j=I-2z \, m_j (0,0,1)^t (1,0,0)$
  in the definition \eref{eq:k-step-matrix} of $a^{(2k+1)}$,
  and collect the lowest and highest powers of~$z$.
\end{proof}

\begin{corollary}
  With $[\cdot]$ denoting the coefficient of the highest occurring power of $z$,
  \begin{equation}
    \label{eq:mrecover}
    m_{n-k} = -\frac{ \bigl[ a^{(2k+1)}_{31}(z) \bigr]}
                    {2 \bigl[ a^{(2k+1)}_{33}(z) \bigr]}, \quad 0\le k\le n-1,
  \end{equation}
  \begin{equation}
    \label{eq:lrecover}
    l_{n-k}=\frac{2 \bigl[ a^{(2k+1)}_{33}(z) \bigr]}
                 {\bigl[ a^{(2k+1)}_{32}(z) \bigr]}, \quad 1\le k\le n-1.
  \end{equation}
\end{corollary}

Next, we show that given $W$ and $Z$, the properties in
Theorem~\ref{thm:a-approx} determine $a(z)$ uniquely.
Like in the DP case \cite{ls-imrp} we are dealing with a problem
of simultaneous rational approximation, with approximants having a common
denominator and for which the functions to be
approximated are related by the quadratic constraint \eref{eq:ZW-relation}.
In addition to the actual approximation conditions \eref{eq:WZapprox},
there is also the symmetry condition \eref{eq:WZsymmetry} which
makes the approximants unique despite the low order of the approximation condition.
This type of approximation seems to be quite distinct
from the cases of classical Hermite--Pad\'e approximation known to us.

Before we proceed to the main body of the paper
we need a few preparatory lemmas concerning
functions $f(z)$ in the complex plane having a Laurent expansion around
$z=\infty$ with finite principal part,
\begin{equation*}
  f(z)=\sum_{j=-\infty}^N f_j z^j.
\end{equation*}
For such functions we let $\pipos$ and $\pineg$ denote the projection operators
onto the subspaces of nonnegative and negative powers, respectively.
For a fixed $\lambda$ we denote by $\frac{1}{z-\lambda}$ the Laurent series
$z^{-1}+\lambda z^{-2}+\lambda^2 z^{-3}+\cdots$.
Finally, for $0\le l$, let $\pil$ denote the projection onto $z^{-l}$.
The following two lemmas are simple exercises in the algebra of Laurent series.

\begin{lemma}
  \label{lem:projections1}
  Suppose $0\le j$ and $1\le l$. Then
  \begin{equation*}
    \pipos \frac{z^j}{z-\lambda}=\frac{z^j-\lambda^j}{z-\lambda},
  \end{equation*}
  \begin{equation*}
    \pineg\frac{z^j}{z-\lambda}=\frac{\lambda^j}{z-\lambda},
    \qquad
    \pil\frac{z^j}{z-\lambda}=\frac{\lambda^{j+l-1}}{z^l}.
  \end{equation*}
  In particular, if $p(z)$ is a polynomial then
  \begin{equation*}
    \pipos \frac{p(z)}{z-\lambda}=\frac{p(z)-p(\lambda)}{z-\lambda},
    \qquad
    \pil \frac{p(z)}{z-\lambda}=\frac{\lambda^{l-1}p(\lambda)}{z^l}.
  \end{equation*}
  Moreover,
  \begin{equation*}
    \Pi ^0  \frac{p(z)}{z-\lambda}=\frac{p(\lambda)-p(0)}{\lambda}.
  \end{equation*}
\end{lemma}

\begin{lemma}
  \label{lem:projections2}
  Let $\kappa\ne \lambda$ be given. Then
  \begin{equation*}
    \pil\frac{z^j}{(z-\kappa)(z-\lambda)}=\frac{\kappa^{j+l-1}-
      \lambda^{j+l-1}}{\kappa-\lambda} z^{-l}.
  \end{equation*}
  In particular, if $p(z)$ is a polynomial then
  \begin{equation*}
    \pil\frac{p(z)}{(z-\kappa)(z-\lambda)}=\frac{\kappa^{l-1}p(\kappa)-
      \lambda^{l-1}p(\lambda)}{\kappa-\lambda} z^{-l}.
  \end{equation*}
  Moreover,
  \begin{equation*}
    \Pi ^0\frac{p(z)}{(z-\kappa)(z-\lambda)}=\frac{1}{\kappa -\lambda}
    \left [ \frac{p(\kappa)-p(0)}{\kappa}-\frac{p(\lambda)-p(0)}{\lambda} \right ]
  \end{equation*}
\end{lemma}

We will introduce three types of approximation problem, each type
specifically referring to the pertinent column of $a(z)$ as in
Theorem~\ref{thm:a-approx}.
The following setup is common to all three cases.

\begin{definition}
  \label{def:setup}
  Let $\{ \lambda_k,b_k,M \}$ be arbitrary numbers satisfying the constraints
  \begin{equation*}
    0 < \lambda_1 < \lambda_2 < \cdots < \lambda_{n-1},
    \qquad
    b_1,\ldots,b_{n-1}<0,
    \qquad
    M>0.
  \end{equation*}
  Define (signed) measures
  \begin{equation}\label{eq:mu-nu}
    \fl
    \rmd\mu(\lambda) = \sum_{k=1}^{n-1} b_k \, \delta(\lambda-\lambda_k), \qquad
    \rmd\nu(\lambda) =-\frac{1}{2M} \delta (\lambda)+ \sum_{k=1}^{n-1} c_k \, \delta(\lambda-\lambda_k),
  \end{equation}
  where the $c_k$'s are given by \eref{eq:cb-relation}, which amounts
  to
  \begin{equation}
    \label{eq:mu-nu2}
    -\int f(\lambda)\,\rmd\nu(\lambda) =
    \frac{f(0)}{2M}
    + \iint \frac{f(\lambda)}{\kappa+\lambda} \rmd\mu(\kappa) \rmd\mu(\lambda)
  \end{equation}
  for any function~$f$. Define the two functions
  \begin{equation}
    \label{eq:WZparfrac}
    \eqalign{
      W(z) = \int \frac{\rmd\mu(\lambda)}{z-\lambda}
      = \sum_{k=1}^{n-1} \frac{b_k}{z-\lambda_k},
      \\
      Z(z) = \int \frac{\rmd\nu(\lambda)}{z-\lambda}
      = -\frac{1/2M}{z} + \sum_{k=1}^{n-1} \frac{c_k}{z-\lambda_k},
    }
  \end{equation}
  and let
  \begin{equation}
    \label{eq:WZstar}
    \fl
    W^*(z) = -W(-z) = \int \frac{\rmd\mu(\lambda)}{z+\lambda},
    \qquad
    Z^*(z) = Z(-z) = -\int \frac{\rmd\nu(\lambda)}{z+\lambda}.
  \end{equation}
  (Note that $Z(z)$ is determined by $W(z)$ and $M$,
  and that $W$ and $Z$ satisfy \eref{eq:ZW-relation}.)
  Finally, let
  \begin{equation}
    \label{eq:Ibeta}
    \beta_j=\int \lambda^j \, \rmd\mu(\lambda),
    \qquad
    I_{ij} = I_{ji} = \iint \frac{\kappa^i \, \lambda^j}{\kappa+\lambda} \rmd\mu(\kappa) \rmd\mu(\lambda).
  \end{equation}
\end{definition}

\begin{definition}
  \label{def:dets}
  Let ${\cal A}_0={\cal B}_0={\cal C}_0={\cal D}_0=1$,
  ${\cal A}_1=I_{00}+\frac{1}{2M}$, ${\cal D}'_1=\beta_0$, and
  \begin{equation}
    {\cal A}_k=\begin{vmatrix}
      I_{00}+\frac{1}{2M} & I_{01} & \cdots & I_{0,k-1}\\
      I_{10} & I_{11} & \cdots & I_{1,k-1} \\
      I_{20} & I_{21} & \cdots & I_{2,k-1} \\
      \vdots & \vdots && \vdots \\
      I_{k-1,0} & I_{k-1,1} & \cdots & I_{k-1,k-1}
    \end{vmatrix}, \qquad 2\le k,
  \end{equation}
  \begin{equation}
    {\cal B}_k = \begin{vmatrix}
      I_{00} & I_{01} & \cdots & I_{0,k-1}\\
      I_{10} & I_{11} & \cdots & I_{1,k-1} \\
      \vdots & \vdots && \vdots \\
      I_{k-1,0} & I_{k-1,1} & \cdots & I_{k-1,k-1}
    \end{vmatrix}, \qquad 1\le k,
  \end{equation}
  \begin{equation}
    {\cal C}_k=\begin{vmatrix}
      I_{11}& I_{12} & \cdots & I_{1k} \\
      I_{21} & I_{22} & \cdots & I_{2k} \\
      \vdots & \vdots && \vdots \\
      I_{k1} & I_{k2} & \cdots & I_{kk}
    \end{vmatrix}, \qquad 1\le k,
  \end{equation}
  \begin{equation}
    {\cal D}_k=\begin{vmatrix}
      I_{10} & I_{11} & \cdots & I_{1,k-1} \\
      I_{20} & I_{21} & \cdots & I_{2,k-1} \\
      \vdots & \vdots && \vdots \\
      I_{k0}& I_{k1} & \cdots & I_{k,k-1}
    \end{vmatrix}, \qquad 1\le k,
  \end{equation}
  \begin{equation}
    {\cal D}'_k=\begin{vmatrix}
      \beta_0& I_{10} & \cdots & I_{1,k-2} \\
      \beta_1 & I_{20} & \cdots & I_{2,k-2} \\
      \vdots & \vdots && \vdots \\
      \beta_{k-1} & I_{k0} & \cdots & I_{k,k-2}
    \end{vmatrix}, \qquad 2\le k.
  \end{equation}
  In all these cases, the index $k$ agrees with the size $k\times k$ of the determinant.
  Note also that ${\cal A}_k = {\cal B}_k + \frac{1}{2M} {\cal C}_{k-1}$ for $k\ge 1$.
\end{definition}

\begin{theorem}[Type~III approximation problem]
  \label{thm:approx-problem3}
  Fix $1\le k\le n-1$.
  There are unique polynomials $Q(z)$, $P(z)$, $\widehat{P}(z)$
  of degree $k$, $k-1$, $k-1$, respectively, satisfying the following
  properties:
  \begin{enumerate}
  \item Approximation:
    \begin{equation}
      \label{eq:PQapprox}
      Q\begin{pmatrix}Z\\W \end{pmatrix}- \begin{pmatrix}\widehat{P}\\
        P\end{pmatrix}= \Or(z^{-1}),
      \quad\text{as $z\to\infty$}.
    \end{equation}
  \item Symmetry:
    \begin{equation}
      \label{eq:PQsymmetry}
      \widehat{P} + P W^* + Q Z^* = \Or(z^{-(k+1)}),
      \quad\text{as $z\to\infty$}.
    \end{equation}
  \item Normalization: $Q(0)=1$.
  \end{enumerate}
  The coefficients of $Q(z)=1+\sum_{j=1}^k q_j z^j$ are given by the
  $k\times k$ linear system
  \begin{equation}
    \label{eq:system-for-Q-explicit}
    \fl
    \begin{pmatrix}
      I_{01} & I_{02} & \cdots & I_{0k} \\
      I_{11} & I_{12} & \cdots & I_{1k} \\
      \vdots &&& \vdots \\
      I_{k-1,1} & I_{k-1,2} & \cdots & I_{k-1,k}
    \end{pmatrix}
    \begin{pmatrix}
      q_1 \\ q_2 \\ \vdots \\ q_k
    \end{pmatrix}
    = -
    \begin{pmatrix}
      I_{00}+\frac{1}{2M} \\ I_{10} \\ \vdots \\ I_{k-1,0}
    \end{pmatrix}.
  \end{equation}
  The polynomials $P(z)$ and $\widehat{P}(z)$ are determined from $Q(z)$ through
  \eref{eq:PPhat} below.
\end{theorem}

\begin{proof}
 Upon applying  $\pipos$ and Lemma \ref{lem:projections1} to \eref{eq:PQapprox} we obtain
  \begin{equation}
    \label{eq:PPhat}
    \eqalign{
      P(z) = \pipos (WQ) &= \int \frac{Q(z)-Q(\lambda)}{z-\lambda}\rmd\mu(\lambda),
      \\
      \widehat{P}(z) = \pipos (ZQ) &= \int\frac{Q(z)-Q(\lambda)}{z-\lambda}\rmd\nu(\lambda),
    }
  \end{equation}
  so that $P$ and $\widehat{P}$ are uniquely determined once we know $Q$.
  Inserting this into \eref{eq:PQsymmetry} yields
  \begin{equation}
    \label{eq:holycow}
    \pipos(ZQ) + \pipos(WQ) \, W^* + Q Z^* = \Or(z^{-(k+1)}).
  \end{equation}
  Applying $\pipos$ to the left-hand side yields $\pipos\bigl((Z+Z^*+WW^*)Q\bigr)$,
  which is zero by \eref{eq:ZW-relation}.
  Thus, the nonnegative powers in \eref{eq:holycow} vanish automatically,
  while the vanishing of the powers $z^{-1},\ldots,z^{-k}$ is the
  condition that will determine the $k$ unknown coefficients in $Q(z)$.
  Define the Hankel and Toeplitz operators
  \begin{equation}
    \label{eq:HankelToeplitz}
    H_f = \pineg \circ M_f \circ \pipos,
    \qquad
    T_f = \pipos \circ M_f \circ \pipos,
  \end{equation}
  where $M_f$ denotes multiplication by the function~$f$.
  Applying $\pineg$ to the left-hand side of \eref{eq:holycow} yields
  $H_{Z^*}Q + H_{W^*} T_W Q$,
  so that $Q$ is determined by the linear system
  \begin{equation}
    \label{eq:system-for-Q}
    (H_{Z^*} + H_{W^*} T_W ) Q = \Or(z^{-(k+1)}).
  \end{equation}
  If we restrict the operators to the finite-dimensional subspaces
  spanned by $\{ 1,z,\ldots,z^k \}$ and $\{ z^{-1},\ldots,z^{-k} \}$
  we obtain $k$ equations
  \begin{equation}
    \Pi ^l (H_{Z^*} + H_{W^*} T_W ) Q=0,  \qquad 1\le l\le k.
  \end{equation}
  With the help of Lemma \ref{lem:projections1} this can be rewritten as
  \begin{equation}
    \fl
    \Pi^l \biggl(
    -\int \frac{Q(z)}{z+\lambda}\rmd\nu(\lambda)
    +\iint \frac{Q(z)-Q(\kappa)}{(z-\kappa)(z+\lambda)}\rmd\mu(\kappa)\rmd\mu(\lambda)
    \biggr)=0,
  \end{equation}
  which according to the relation \eref{eq:mu-nu2} between the measures
  $\rmd\mu$ and $\rmd\nu$ equals
  \begin{equation}
    \fl
    \Pi^l \biggl(
    \frac{Q(z)/z}{2M}
    +\iint \biggl(
    \frac{Q(z)/(z+\lambda)}{\kappa+\lambda}
    +\frac{Q(z)-Q(\kappa)}{(z-\kappa)(z+\lambda)}
    \biggr)
    \rmd\mu(\kappa) \rmd\mu(\lambda)
    \biggr)=0.
  \end{equation}
  Finally,
  upon executing the projection $\Pi ^l$ for $1\le l \le k$
  using formulas from Lemmas~\ref{lem:projections1}
  and~\ref{lem:projections2}, we obtain
  \begin{equation}
    \label{eq:ortho}
    \eqalign{
      \frac{Q(0)}{2M}
      +\iint \frac{Q(\kappa)} {\kappa +\lambda}
      \rmd\mu(\kappa) \rmd\mu (\lambda)=0,
      \\
      \iint \frac{(\lambda)^{l-1} Q(\kappa)}
      {\kappa +\lambda} \rmd\mu(\kappa) \rmd\mu (\lambda)=0, \qquad 2\le l\le k.
    }
  \end{equation}
  When $Q$ is represented by the column vector $(1,q_1,\ldots,q_k)^t$ this
  takes the form \eref{eq:system-for-Q-explicit}.
  The determinant of the system \eref{eq:system-for-Q-explicit}
  is ${\cal D}_k$ of Definition~\ref{def:dets},
  which is positive by Lemma~\ref{lem:super-Heine} in Appendix~B.
  Thus the solution is unique.
\end{proof}

\begin{definition}
  \label{def:denote3k}
  Denote the results of the $k$th approximation problem of
  type~III by $Q_{3k}, P_{3k}, \widehat {P}_{3k}$.
\end{definition}

\begin{corollary}
  \label{cor:q3}
  The highest order coefficient of $Q_{3k}(z)$ is
  $[Q_{3k}]=(-1)^k {\cal A}_k/{\cal D}_k$.
\end{corollary}

\begin{theorem}[Type~II approximation problem]
  \label{thm:approx-problem2}
  Fix $1\le k\le n-1$.
  There are unique polynomials $Q(z)$, $P(z)$, $\widehat{P}(z)$
  of degree $k$, $k-1$, $k-1$, respectively, satisfying the following
  properties:
  \begin{enumerate}
  \item Approximation:
    \begin{equation}
      \label{eq:PQapprox2}
      Q\begin{pmatrix}Z\\W \end{pmatrix}- \begin{pmatrix}\widehat{P}\\
        P\end{pmatrix}= \begin{pmatrix}\Or(z^{-1})\\\Or(1)\end{pmatrix},
      \quad\text{as $z\to\infty$}.
    \end{equation}
  \item Symmetry:
    \begin{equation}
      \label{eq:PQsymmetry2}
      \widehat{P} + P W^* + Q Z^* = \Or(z^{-(k+1)}),
      \quad\text{as $z\to\infty$}.
    \end{equation}
  \item Normalization: $Q(0)=0, \,P(0)=1$.    \end{enumerate}
  The coefficients of $Q(z)=\sum_{j=1}^k q_j z^j$ are given by the
  $k\times k$ linear system
  \begin{equation}
    \label{eq:system-for-Q-explicit2}
    \begin{pmatrix}
      I_{10} & I_{11} & \cdots & I_{1,k-1} \\
      I_{20} & I_{21} & \cdots & I_{2,k-1} \\
      \vdots &&& \vdots \\
      I_{k0} & I_{k1} & \cdots & I_{k,k-1}
    \end{pmatrix}
    \begin{pmatrix}
      q_1 \\ q_2 \\ \vdots \\ q_k
    \end{pmatrix}
    =
    \begin{pmatrix}
      \beta_0 \\ \beta_1 \\\vdots \\ \beta_{k-1}
    \end{pmatrix}.
  \end{equation}
  The polynomials $P(z)$ and $\widehat{P}(z)$ are determined from $Q(z)$ through
  \eref{eq:PPhat2} below.
\end{theorem}

\begin{proof}
  The proof proceeds in a similar way to the proof of Theorem~\ref{thm:approx-problem3}.
  We therefore only indicate the steps which are new to the present case.
  First, from the normalization condition we get
  \begin{equation}
    \label{eq:PPhat2}
    \eqalign{
      P(z) = 1+\Pi_{>0} (WQ) = \pipos (WQ)+1-\Pi^0(WQ),
      \\
      \widehat{P}(z) = \pipos (ZQ) = \int\frac{Q(z)-Q(\lambda)}{z-\lambda}\rmd\nu(\lambda).
    }
  \end{equation}
  Inserting this into \eref{eq:PQsymmetry2} yields
  \begin{equation}
    \label{eq:holycow2}
    \fl
    \pipos(ZQ) + \pipos(WQ) \, W^* + Q Z^* + \bigl(1-\Pi^0(WQ)\bigr)\, W^*= \Or(z^{-(k+1)}).
  \end{equation}
  Applying $\pipos$ to the left-hand side yields $\pipos\bigl((Z+Z^*+WW^*)Q\bigr)$,
  which is zero by \eref{eq:ZW-relation}.
  Thus, the nonnegative powers in \eref{eq:holycow2} vanish automatically,
  while the vanishing of the powers $z^{-1},\ldots,z^{-k}$ is the
  condition that will determine the $k$ unknown coefficients in $Q(z)$.
  With the help of the Hankel and Toeplitz operators defined in the previous
  proof we find that $Q$ is determined by the linear system
  \begin{equation}
    \label{eq:system-for-Q2}
    (H_{Z^*} + H_{W^*} T_W ) Q + \bigl(1-\Pi^0(WQ)\bigr)\, W^*= \Or(z^{-(k+1)}),
  \end{equation}
  where we have used that $W^*=\Or(1/z)$ according to its definition.
  If we restrict the operators to the finite-dimensional subspaces
  spanned by $\{ 1,z,\ldots,z^k \}$ and $\{ z^{-1},\ldots,z^{-k} \}$
  we obtain the $k$ equations
  \begin{equation}
    \fl
    \Pi ^l \biggl(
    (H_{Z^*} + H_{W^*} T_W ) Q + \bigl(1-\Pi^0(WQ)\bigr)\, W^*
    \biggr)=0,
    \qquad 1\le l\le k.
  \end{equation}
  With the help of
  Lemma~\ref{lem:projections2} and the relation \eref{eq:mu-nu2} between
  the measures $\rmd\mu$ and $\rmd\nu$ we obtain
  \begin{equation}\label{eq:ortho2}
    \fl
    \iint \frac{\lambda ^l \, Q(\kappa)}
    {\kappa(\kappa +\lambda)} \rmd\mu(\kappa) \rmd\mu (\lambda)
    =\int(\lambda)^{l-1} \rmd\mu(\lambda),
    \qquad 1\le l\le k.
  \end{equation}
  When $Q$ is represented by the column vector $(q_1,\ldots,q_k)^t$ this
  gives the system \eref{eq:system-for-Q-explicit2}, which
  has the same determinant ${\cal D}_k>0$ as \eref{eq:system-for-Q-explicit}.
  Hence, the solution is again unique.
\end{proof}

\begin{definition}
  \label{def:denote3kplus1}
  Denote the results of the $k$th approximation problem of type~II by
  $Q_{3k+1}, P_{3k+1}, \widehat {P}_{3k+1}$.
\end{definition}

\begin{corollary}\label{cor:q2}
  The highest order coefficient of $Q_{3k+1}(z)$ is
  $[Q_{3k+1}]=(-1)^{k-1} {\cal D}'_k/{\cal D}_k$.
\end{corollary}

\begin{theorem}[Type~I approximation problem]
  \label{thm:approx-problem1}
  Let $0\le k\le n-1$.
  Then there are unique polynomials $Q(z)$, $P(z)$, $\widehat{P}(z)$
  of degree $k+1$, $k$, $k$, respectively, satisfying the following
  properties:
  \begin{enumerate}
  \item Approximation:
    \begin{equation}
      \label{eq:PQapprox1}
      Q\begin{pmatrix}Z\\W \end{pmatrix}- \begin{pmatrix}\widehat{P}\\
        P\end{pmatrix}= \begin{pmatrix}\Or(z^{-1})\\\Or(1)\end{pmatrix},
      \quad\text{as $z\to\infty$}.
    \end{equation}
  \item Symmetry:
    \begin{equation}
      \label{eq:PQsymmetry1}
      \widehat{P} + P W^* + Q Z^* = \Or(z^{-(k+1)}),
      \quad\text{as $z\to\infty$}.
    \end{equation}
  \item Normalization: $Q(0)=0, \,P(0)=0, \, \widehat{P}(0)=1$.    \end{enumerate}
  The coefficients of $Q(z)=\sum_{j=1}^{k+1} q_j z^j$ are given by the
  $(k+1) \times (k+1)$ linear system
  \begin{equation}
    \label{eq:system-for-Q-explicit1}
    \begin{pmatrix}
      I_{00}+\frac{1}{2M} & I_{01} & \cdots & I_{0k}\\
      I_{10} & I_{11} & \cdots & I_{1k} \\
      I_{20} & I_{21} & \cdots & I_{2k} \\
      \vdots &&& \vdots \\
      I_{k0} & I_{k1} & \cdots & I_{kk}
    \end{pmatrix}
    \begin{pmatrix}
      q_1 \\ q_2 \\q_3\\ \vdots \\ q_{k+1}
    \end{pmatrix}
    =
    -\begin{pmatrix}
      1\\ 0\\0\\\vdots \\ 0
    \end{pmatrix}.
  \end{equation}
  The polynomials $P(z)$ and $\widehat{P}(z)$ are determined from $Q(z)$ through
  \eref{eq:PPhat1} below.
\end{theorem}

\begin{proof}
  From the normalization and approximation conditions we get
  \begin{equation}
    \label{eq:PPhat1}
    \eqalign{
      P(z) = \Pi_{>0} (WQ)=&\pipos (WQ)-\Pi^0(WQ),
      \\
      \widehat{P}(z) = \pipos (ZQ)=&\int\frac{Q(z)-Q(\lambda)}{z-\lambda}\rmd\nu(\lambda).
    }
  \end{equation}
  Inserting this into \eref{eq:PQsymmetry1} yields
  \begin{equation}
    \label{eq:holycow1}
    \fl
    \pipos(ZQ) + \pipos(WQ) \, W^* + Q Z^* -\Pi^0(WQ)\, W^*= \Or(z^{-(k+1)}).
  \end{equation}
  Applying $\pipos$ to the left-hand side yields $\pipos\bigl((Z+Z^*+WW^*)Q\bigr)$,
  which is zero by \eref{eq:ZW-relation}.
  Thus, the nonnegative powers in \eref{eq:holycow1} vanish automatically,
  while the vanishing of the powers $z^{-1},\ldots,z^{-k}$ will give $k$ conditions
  on the $k+1$ unknown coefficients in $Q(z)$, which together with the normalization
  condition $\widehat{P}(0)=1$ is the required number of conditions.
  Symbolically, with the help of the Hankel and Toeplitz operators defined
  earlier, we find that $Q$ is determined by the linear system
  \begin{equation}
    \label{eq:system-for-Q1}
    \eqalign{
      (H_{Z^*} + H_{W^*} T_W ) Q -\Pi^0(WQ)\, W^*= \Or(z^{-(k+1)}), \\
      \Pi ^0(ZQ)=1.
    }
  \end{equation}
  If we restrict the operators to the finite-dimensional subspaces
  spanned by $\{z,\ldots,z^{k+1} \}$ and $\{ z^{-1},\ldots,z^{-k} \}$,
  then the $k$ equations alluded to earlier are
  \begin{equation}
    \fl
    \Pi ^l \biggl(
    (H_{Z^*} + H_{W^*} T_W ) Q-\Pi^0(WQ)\, W^*
    \biggr)=0,
    \qquad 1\le l\le k.
  \end{equation}
  Likewise, Lemma~\ref{lem:projections2} and the relation \eref{eq:mu-nu2} between
  the measures $\rmd\mu$ and $\rmd\nu$ imply that the complete set of equations
  determining~$Q$ reads
  \begin{equation}
    \label{eq:ortho1}
    \eqalign{
      \int \frac{Q(\lambda)}{\lambda} \rmd\nu(\lambda)=1,
      \\
      \iint \frac{\lambda ^l \, Q(\kappa)}
      {\kappa(\kappa +\lambda)} \rmd\mu(\kappa) \rmd\mu (\lambda)= 0,
      \qquad 1\le l\le k.
    }
  \end{equation}
  When $Q$ is represented by the column vector $(q_1,\ldots,q_{k+1})^t$ this
  gives \eref{eq:system-for-Q-explicit1}.  The determinant of the
  matrix of coefficients is
  ${\cal A}_{k+1} = {\cal B}_{k+1} + \frac{1}{2M} {\cal C}_k$
  in the notation of Definition~\ref{def:dets}.
  Following the methods used in
  \cite[Lemma 4.10]{ls-imrp} one can establish Heine-type formulas for both
  determinants in this sum, which show that both are strictly positive provided
  $k<n-1$; see Lemmas~\ref{lem:super-Heine2} and ~\ref{lem:super-Heine3} in Appendix~B.
  From these formulas one also concludes that
  when $k=n-1$ the first determinant is zero, since the measure
  $\rmd\mu(\lambda)$ is supported on $n-1$ points, while the second remains
  strictly positive.
  This proves that the system \eref{eq:system-for-Q-explicit1} has a
  unique solution for the claimed range of $k$.
\end{proof}

\begin{definition}
  \label{def:denote3kplus2}
  Denote the results of the $k$th approximation problem of type~I by
  $Q_{3k+2}, P_{3k+2}, \widehat {P}_{3k+2}$.
\end{definition}

\begin{corollary}
  \label{cor:q1}
  The highest order coefficient of $Q_{3k+2}(z)$ is
  $[Q_{3k+2}]=(-1)^{k+1} {\cal C}_k/{\cal A}_k$.
\end{corollary}

The last approximation problem of type~I has a particular significance.

\begin{corollary}
  \label{cor:laststep}
  For $k=n-1$, the approximation problem of type~I gives
  \begin{equation}
    Q_{3n-1}(z)=-2Mz \prod_{j=1}^{n-1} \left(1-\frac{z}{\lambda_j}\right).
  \end{equation}
  Moreover, the approximation is exact:
  \begin{equation}
    Z(z)=\frac{\widehat{P}_{3n-1}(z)}{Q_{3n-1}(z)}, \qquad W(z)=\frac{P_{3n-1}(z)}{Q_{3n-1}(z)}.
  \end{equation}
\end{corollary}

\begin{proof}
  We write the conditions \eref{eq:ortho1} that $Q_{3n-1}$ has to satisfy:
  \begin{equation}
    \label{eq:ortho1again}
    \eqalign{
      \int \frac{Q_{3n-1}(\lambda)}{\lambda} \rmd\nu(\lambda)=1,\\
      \iint \frac{\lambda^i \, Q_{3n-1}(\kappa)}
      {\kappa(\kappa +\lambda)} \rmd\mu(\kappa) \rmd\mu (\lambda)= 0,
      \qquad 1\le i\le n-1.
    }
  \end{equation}
  We claim that $Q_{3n-1}$ has to vanish on the support of the
  measure $\rmd\mu$. We recall that $Q_{3n-1}$ is a polynomial of
  degree $n$ with one root at $\lambda =0$ by the definition of
  the approximation problem of type~I. Thus the claim is that the
  remaining roots are precisely the $\lambda_j$ occuring in \eref{eq:mu-nu}.
  To prove that, we write explicitly the second condition,
  \begin{equation*}
    \sum_{j=1}^{n-1}
    \left( \int \frac{\lambda^i}{ \lambda +\lambda_j}\rmd\mu(\lambda) \right)
    \left( \frac{Q_{3n-1}(\lambda_j) \, b_j}{\lambda_j} \right) =0,
    \qquad 1\le i\le n-1.
  \end{equation*}
  It is straightforward to establish a multiple integral formula of
  Heine type for the matrix of coefficients, from which one concludes
  that the matrix is nonsingular;
  see \eref{eq:Cauchy-type} in Appendix~B.
  This implies that $Q_{3n-1}(\lambda_j)=0$ for $1\le j\le n-1$, which proves the claim.
  Moreover, the definition of the measure $\rmd\nu$ and the first
  condition in \eref{eq:ortho1again}
  imply that the coefficient of the linear term of $Q_{3n-1}$
  is $-2M$, which proves the product formula for $Q_{3n-1}$.
  To prove that the $k=n-1$ approximation of type~I is exact we show that
  \begin{equation*}
    \pineg \left(
      Q_{3n-1} \begin{pmatrix} Z\\W \end{pmatrix}
      - \begin{pmatrix} \widehat{P}_{3n-1} \\ P_{3n-1} \end{pmatrix}
    \right)
    = 0.
  \end{equation*}
  This follows from Lemma~\ref{lem:projections1}, since for example
  \begin{equation*}
    \pineg (Q_{3n-1}W)(z)
    =\int\frac{Q_{3n-1}(z)}{z-\lambda}\rmd\mu(\lambda)
    =\int \frac{Q_{3n-1}(\lambda)}{z-\lambda}\rmd\mu (\lambda) =0
  \end{equation*}
  because $Q_{3n-1}$ vanishes on the support of $\rmd\mu$.  Likewise,
  \begin{equation*}
    \Pi ^0 (Q_{3n-1}W)
    =\Pi ^0 \int \frac{Q_{3n-1}(z)}{z-\lambda} \rmd\mu(\lambda)
    =\int \frac{Q_{3n-1}(\lambda)}{\lambda} \rmd\mu(\lambda)=0.
  \end{equation*}
  The rest follows from the definition of
  the approximation problem of type~I.
\end{proof}

We can now present the solution of the uniqueness part of the inverse problem.

\begin{theorem}
  \label{thm:recovery}
  Let
  \begin{equation}
    \label{eq:constraints}
    0 < \lambda_1 < \lambda_2 < \cdots < \lambda_{n-1},
    \qquad
    b_1,\ldots,b_{n-1}<0,
    \qquad
    M>0
  \end{equation}
  be the spectral data (as described in Section~\ref{sec:spectral})
  of a cubic string $\{m_k,x_k\}_{k=1}^n$ with $m_k>0$ and $x_1<\cdots<x_n$.
  Then the string data can be recovered uniquely from the spectral data,
  up to a translation along the $x$ axis, from the determinantal formulas
  \begin{equation}
    \label{eq:l-m}
    m_{n-k} = \frac{{\cal C}_k {\cal D}_k}{2{\cal A}_{k+1}{\cal A}_k},
    \qquad
    l_{n-k} = -\frac{2 {\cal A}_k}{{\cal D}'_k},
  \end{equation}
  where $l_j=x_{j+1}-x_j$.
  (See Definitions~\ref{def:setup} and~\ref{def:dets} for notation.)
\end{theorem}

\begin{proof}
  Equations \eref{eq:l-m} follow from
  Corollaries \ref{cor:q3}, \ref{cor:q2}, \ref{cor:q1} and
  equations \eref{eq:mrecover} and \eref{eq:lrecover}.
\end{proof}

To show the existence part (that any numbers satisfying the constraints
\eref{eq:constraints} really are spectral data of some cubic string),
we need the following theorem, whose proof will be published elsewhere.
Note that the orthogonal polynomials occurring in the description of continued
fractions of Stieltjes type satisfy three-term recurrence relations.
It is therefore gratifying to know that the polynomials occurring in
the approximation problems for the cubic string satisfy four-term
recurrence relations.

\begin{theorem}
  \label{thm:4recurrence}
  Consider the same setup as in Definition~\ref{def:setup}
  (numbers $\{ \lambda_k,b_k,M \}$ given),
  and let the polynomials $Q_j(z)$, $P_j(z)$, $\widehat{P}_j(z)$ be defined (for $j\ge 2$)
  by the unique solutions of the approximation problems
  of type~I, II and~III,
  as in Definitions \ref{def:denote3k}, \ref{def:denote3kplus1} and~\ref{def:denote3kplus2}.
  Define $l_j$ and $m_j$ from the numbers $\{ \lambda_k, b_k, M \}$ via \eref{eq:l-m},
  and consider the recurrence relations
  \begin{equation}
    \label{eq:4recurrence}
    \fl
    \eqalign{
      X_{3k\phantom{+1}} = \textstyle\frac12 l_{n-k}^2 \, X_{3k-1} + l_{n-k}\, X_{3k-2} + X_{3k-3},\qquad & 1\le k\le n-1,
      \\
      X_{3k+1} = l_{n-k}\, X_{3k-1} + X_{3k-2}, & 1\le k\le n-1,
      \\
      X_{3k+2} = -2z\, m_{n-k}\, X_{3k} + X_{3k-1}, & 0\le k\le n-1.
    }
  \end{equation}
  Then with initial conditions $(X_{-1},X_0,X_1)=(1,0,0)$
  the solution (for $j\ge 2$) is $X_j = \widehat{P}_j$,
  when $(X_{-1},X_0,X_1)=(0,1,0)$
  the solution is $X_j = Q_j$,
  and
  when $(X_{-1},X_0,X_1)=(0,0,1)$
  the solution is $X_j = P_j$.
\end{theorem}

\begin{theorem}
  \label{thm:existence}
  Given any numbers $\{ \lambda_k,b_k,M \}$
  satisfying the constraints \eref{eq:constraints},
  the formulas \eref{eq:l-m} define a cubic string with these
  numbers as spectral data.
\end{theorem}

\begin{proof}
  Note first that all $m_j$ and $l_j$ defined by \eref{eq:l-m}
  are strictly positive.
  This follows since ${\cal D}'_k<0$ according to \eref{eq:super-Heine2} and
  our definition of $\rmd\mu(\lambda)$ as a negative measure,
  while the remaining determinants are strictly positive.
  This means that we have a valid cubic string,
  but we do not know yet whether its spectral data
  $\{ \widetilde{\lambda}_k, \widetilde{b}_k, \widetilde{M} \}$
  coincide with the numbers $\{ \lambda_k,b_k,M \}$ that we started with.
  However, plugging
  $\{ \widetilde{\lambda}_k, \widetilde{b}_k, \widetilde{M} \}$
  into \eref{eq:l-m} recovers the same string data $\{m_j,l_j\}$
  by uniqueness (Theorem~\ref{thm:recovery}).
  Thus, applying Theorem~\ref{thm:4recurrence}
  to the set of numbers $\{ \lambda_k,b_k,M \}$ 
  produces the same polynomials $Q_j$, $P_j$, $\widehat{P}_j$
  as one would get with
  $\{ \widetilde{\lambda}_k, \widetilde{b}_k, \widetilde{M} \}$,
  since the string data are the same in both cases.
  And according to Corollary~\ref{cor:laststep},
  both $M$ and $W(z)$ are uniquely determined from these polynomials.
  Thus 
  $\widetilde{\lambda}_k = \lambda_k$, $\widetilde{b}_k = b_k$,
  and $\widetilde{M} = M$, which is what we needed to show.
\end{proof}

\ack
JK wishes to acknowledge the support received from the National
Science and Engineering Research Council of Canada (NSERC) in the
form of the Undergraduate Student Research Award.
HL is grateful to the Institut Mittag-Leffler, Djursholm, Sweden,
for its hospitality during the Fall 2005 program ``Wave Motion''.
The research of JS is supported by NSERC.

\appendix

\section{A collection of facts about totally nonnegative and oscillatory matrices}
\def\thesection{\Alph{section}} 

We proved that the spectrum of the Neumann-like discrete cubic string
\eref{eq:cubic-neumann} is positive and simple by rewriting it as a
matrix eigenvalue problem with an oscillatory matrix.
Below, the reader will find some relevant definitions and theorems
from the vast area of totally nonnegative matrices;
for details see the book by Gantmacher and Krein \cite{gantmacher-krein}
and the excellent survey by Fomin and Zelevinsky \cite{Fomin}.

\begin{definition}
  A square matrix is \emph{totally nonnegative (TN)} if all its minors
  are nonnegative. It is \emph{totally positive (TP)} if all its
  minors are positive. It is \emph{oscillatory} if it is TN and some
  power of it is TP.
\end{definition}

\begin{theorem}
  \label{thm:TPspectrum}
  All eigenvalues of a TP matrix are positive and of algebraic
  multiplicity one, and likewise for oscillatory matrices. All
  eigenvalues of a TN matrix are nonnegative, but in general of
  arbitrary multiplicity.
\end{theorem}

\begin{theorem}
  \label{thm:product}
  The product of an oscillatory matrix and a nonsingular TN matrix is
  oscillatory.
\end{theorem}

\begin{theorem}[Star operation]
  \label{thm:signoscillatory}
  If $A=(a_{ij})$ is oscillatory, then $(A^*)^{-1}$ is oscillatory,
  where $A^*=((-1)^{i+j}a_{ij})$.
\end{theorem}

\begin{theorem}
  \label{thm:tridiagonal}
  A symmetric tridiagonal matrix is oscillatory if and only if it is
  positive definite and has positive off-diagonal entries.
\end{theorem}

\begin{definition}
  A planar network $(\Gamma, \omega)$ of order $n$ is an acyclic, planar
  directed graph~$\Gamma$ with arrows going from left to right,
  with $n$ sources (vertices with outgoing arrows only) on the left side,
  and with $n$ sinks (vertices with incoming arrows only) on the right side.
  The sources and sinks are numbered $1$ to~$n$ from top to bottom.
  All other vertices have at least one arrow coming in and at least one arrow going out.
  Each edge~$e$ of the graph~$\Gamma$ is assigned a scalar weight~$\omega(e)$.
\end{definition}

\begin{definition}
  \label{weightmatrix}
  Given a planar network $\Gamma$ of order $n$, the
  \emph{weight} of a directed path in $\Gamma$ is the product of all the weights
  of the edges of that path.
  The \emph{weighted path matrix} $A(\Gamma,\omega)$ is
  the $n\times n$ matrix whose $(i,j)$ entry is the sum of
  the weights of the possible paths from source~$i$ to sink~$j$.
\end{definition}

The following famous theorem was discovered independently by
Lindstr\"om, Karlin--McGregor, and Gessel--Viennot.

\begin{theorem}
  \label{lem:Lind}
  Let $I$ and $J$ be subsets of $\{1,\ldots,n\}$ with the same cardinality.
  The minor
  \begin{equation*}
    \Delta_{I,J}(A) = \det (a_{ij})_{i\in I,\,j\in J}
  \end{equation*}
  of the weighted path matrix $A(\Gamma,\omega)$ of a
  planar network is equal to the sum of the weights of
  all possible families of vertex-disjoint paths connecting the sources labelled by~$I$ to
  the sinks labelled by~$J$.
  (The weight of a family of paths is defined as
  the product of the weights of the individual paths.)
\end{theorem}

\begin{corollary}
  \label{TNcorollary}
  If all weights $\omega$ are nonnegative, then the weighted path matrix is TN.
\end{corollary}

\def\thesection{Appendix \Alph{section}} 
\section{Heine-like determinant formulas}
\def\thesection{\Alph{section}} 

The following was proved in \cite{ls-imrp},
where it appeared as Lemma~4.10.

\begin{lemma}
  \label{lem:super-Heine}
  Suppose $\mu$ is a measure on $\R$ such that the integrals
  \begin{equation}
    \label{eq:beta-I}
    \beta_a = \int x^a \, \rmd\mu(x),
    \qquad
    I_{ab} =
    \iint \frac{x^{a} y^{b}}{x+y} \,\rmd\mu(x) \rmd\mu(y)
  \end{equation}
  are finite. Let
  \begin{equation}
    \label{eq:integrals-uv}
    \eqalign{
      u_k &=
      \frac{1}{k!}
      \int_{\R^k} \frac{\Delta (x)^2}{\Gamma(x)} \rmd\mu^k(x),
      \\
      v_k &=
      \frac{1}{k!}
      \int_{\R^k} \frac{\Delta(x)^2}{\Gamma(x)} \, x_1 x_2 \ldots x_k \, \rmd\mu^k(x),
    }
  \end{equation}
  where
  \begin{equation}
    \eqalign{
      \Delta (x)&=\Delta(x_1,\ldots,x_k)=\prod_{i<j}(x_i-x_j),\\
      \Gamma(x)&=\Gamma(x_1,\ldots,x_k)=\prod_{i<j}(x_i+x_j).
    }
  \end{equation}
  (When $k=0$ or~$1$, we let $\Delta(x)=\Gamma(x)=1$.
  Likewise, $u_0=v_0=1$.)
  Then for $k \ge 1$ the following $k\times k$ determinant formulas hold:
  \begin{equation}
    \label{eq:super-Heine1}
    {\cal D}_k :=
    \begin{vmatrix}
      I_{10} & I_{11} & I_{12} & \cdots & I_{1,k-1} \\
      I_{20} & I_{21} & I_{22} & \cdots & I_{2,k-1} \\
      I_{30} & I_{31} & I_{32} & \cdots & I_{3,k-1} \\
      \vdots & \vdots & \vdots && \vdots \\
      I_{k0} & I_{k1} & I_{k2} & \cdots & I_{k,k-1}
    \end{vmatrix}
    = \frac{(u_k)^2}{2^k},
  \end{equation}
  \begin{equation}
    \label{eq:super-Heine2}
    {\cal D}'_k :=
    \begin{vmatrix}
      \beta_0 & I_{10} & I_{11} & \cdots & I_{1,k-2} \\
      \beta_1 & I_{20} & I_{21} & \cdots & I_{2,k-2} \\
      \beta_2 & I_{30} & I_{31} & \cdots & I_{3,k-2} \\
      \vdots & \vdots & \vdots && \vdots \\
      \beta_{k-1} & I_{k0} & I_{k1} & \cdots & I_{k,k-2}
    \end{vmatrix}
    = \frac{u_k \, u_{k-1}}{2^{k-1}},
  \end{equation}
  \begin{equation}
    \label{eq:super-Heine3}
    {\cal D}''_k :=
    \begin{vmatrix}
      \beta_0 & I_{11} & I_{12} & \cdots & I_{1,k-1} \\
      \beta_1 & I_{21} & I_{22} & \cdots & I_{2,k-1} \\
      \beta_2 & I_{31} & I_{32} & \cdots & I_{3,k-1} \\
      \vdots & \vdots & \vdots && \vdots \\
      \beta_{k-1} & I_{k1} & I_{k2} & \cdots & I_{k,k-1}
    \end{vmatrix}
    = \frac{u_k \, v_{k-1}}{2^{k-1}}.
  \end{equation}
\end{lemma}

The following lemmas state some additional determinant formulas
needed in this paper. See Definition~\ref{def:dets} for notation.

\begin{lemma}
  \label{lem:super-Heine2}
  For $k\ge 1$,
  \begin{equation}
    \label{eq:Bk}
    {\cal B}_k
    = \frac{1}{(2k)!}
    \int_{\R^{2k}} 
    \frac{1}{\Gamma(x)} \, 
    \left(
      \sum_{I} \Delta_I^2 \Delta_{I^c}^2 \Gamma_I \Gamma_{I^c}
    \right)
    \, \rmd\mu^{2k}(x),
  \end{equation}
  where
  \begin{equation}
    \Delta_I^2=\Delta(x_{i_1},\dots,x_{i_k})^2,
    \qquad
    \Gamma_I=\Gamma(x_{i_1},\dots,x_{i_k}),
  \end{equation}
  where the sum runs over all $\binom{2k}{k}$ $k$-element subsets $I = \{ i_1,\dots,i_k \}$
  of $\{ 1,\dots,2k \}$,
  and where $I^c$ denotes the complement of~$I$.
  Moreover, letting
  \begin{equation}
    t_k = \frac{1}{k!} \int_{\R^k} \frac{\Delta (x)^2}{\Gamma(x)} \,
    \frac{\rmd\mu^k(x)}{x_1 x_2 \ldots x_k}
    \quad\text{for $k\ge 1$},
  \end{equation}
  we have
  \begin{equation}
    \label{eq:w-like-thing}
    {\cal B}_k = \frac{1}{2^k}
    \begin{vmatrix}
      t_k & u_{k-1} \\
      t_{k+1} & u_k
    \end{vmatrix}.
  \end{equation}
\end{lemma}

The proof of \eref{eq:Bk} is omitted since it is very similar to the proofs in \cite{ls-imrp}.
Equation \eref{eq:w-like-thing} follows by rewriting the integrand in \eref{eq:Bk} using
the symmetric function identity (2.76) from \cite{ls-imrp}.
By simply replacing the measure $\rmd\mu(x)$ by $x\,\rmd\mu(x)$
we obtain the corresponding formulas for ${\cal C}_k$:

\begin{lemma}
  \label{lem:super-Heine3}
  For $k\ge 1$,
  \begin{equation}
    \label{eq:Ck}
    {\cal C}_k
    = \frac{1}{(2k)!}
    \int_{\R^{2k}} 
    \frac{x_1 x_2 \ldots x_{2k}}{\Gamma(x)} \, 
    \left(
      \sum_{I} \Delta_I^2 \Delta_{I^c}^2 \Gamma_I \Gamma_{I^c}
    \right)
    \, \rmd\mu^{2k}(x).
  \end{equation}
  Moreover,
  \begin{equation}
    {\cal C}_k = \frac{1}{2^k}
    \begin{vmatrix}
      u_k & v_{k-1} \\
      u_{k+1} & v_k
    \end{vmatrix}.
  \end{equation}
\end{lemma}

It follows from \eref{eq:Bk} and \eref{eq:Ck} that if
$\rmd\mu$ is a discrete measure supported on $n-1$ positive points,
like in the main text,
then ${\cal B}_k$ and ${\cal C}_k$ are strictly positive for $1 \le k < n$ and zero for $k\ge n$
(because then each $\Delta_I$ vanishes).

Finally we state a formula for the determinant appearing in the proof of
Corollary~\ref{cor:laststep}.
We omit this proof as well; it is also similar to the previous ones
(and actually a little easier, since no double integrals are involved).

\begin{lemma}
  \label{lem:Cauchy-type}
  Let the measure $\rmd\mu$ be the same as in the main text and define the
  $(n-1)\times (n-1)$ determinant ${\cal E}=\det[{\cal E}_{ij}]$ where 
  \begin{equation*}
    {\cal E}_{ij}=\int \frac{x^i}{x +\lambda_j}\rmd\mu(x),
    \qquad
    i,\,j=1,\ldots,n-1.
  \end{equation*}
  Then 
  \begin{equation}
    \label{eq:Cauchy-type}
    {\cal E}=
    \frac{\Delta(\lambda_{n-1}, \lambda_{n-2}, \dots, \lambda_1)}{(n-1)!}
    \int_{\R^{n-1}}
    \frac{\Delta(x)^2}{\prod _{i,j=1}^{n-1} (x_i+\lambda_j)} \rmd\mu^{n-1}(x).  
  \end{equation}
  In particular, ${\cal E} \ne 0$.
\end{lemma}

\section*{References}
\bibliographystyle{iopart-num-custom}
\bibliography{cubic_string_IP-ver2}

\end{document}